\documentclass[11pt]{article}

\usepackage{amscd}

\title{Preprojective cluster variables of acyclic cluster algebras\thanks{Supported in part by the NSF of China (Grants 10471071).}}

\author{Bin Zhu\thanks{E-mail: bzhu@math.tsinghua.edu.cn}}

\date{ \small {Department of Mathematical Sciences, \\  Tsinghua University,
100084 Beijing, P. R. China\\
 \vspace{0.4cm} \small{Dedicated to Professor Yingbo Zhang on  the occasion of her 60th birthday}}}

\begin{document}

\maketitle

\def\s{\stackrel}
\def\gama{\gamma}
\def\Longrightarrow{{\longrightarrow}}
\def\P{{\cal P}}
\def\A{{\cal A}}
\def\F{\mathcal{F}}
\def\X{\mathcal{X}}
\def\T{\mathcal{T}}
\def\m{\textbf{ M}}
\def\t{{\tau }}
\def\b{\textbf{d}}
\def\K{{\cal K}}

\def\G{{\Gamma}}
\def\e{\mbox{exp}}

\def\righta{\rightarrow}

\def\s{\stackrel}

\def\ncong{\not\cong}

\def\mathbb{\NN}

\def\Hom{\mbox{Hom}}
\def\Ext{\mbox{Ext}}
\def\ind{\mbox{ind}}
\def\coprod{\amalg }
\def\L{\Lambda}
\def\c{\circ}

\newcommand{\uHom}{\operatorname{\underline{Hom}}\nolimits}
\newcommand{\End}{\operatorname{End}\nolimits}
\renewcommand{\r}{\operatorname{\underline{r}}\nolimits}
\def \text{\mbox}


\begin{center}

\begin{minipage}{12cm}{\footnotesize\textbf{Abstract.} It is proved that any cluster-tilted algebra defined in the cluster category $\mathcal{C}(H)$ has the same representation type as the initial hereditary algebra $H$.
For any valued quiver $(\G,\Omega)$, an injection from the subset $\mathcal{PI}(\Omega)$ of
  the cluster category $\mathcal{C}(\Omega)$ consisting of indecomposable preprojective objects, preinjective objects and the first shifts of indecomposable
 projective modules to the set of cluster variables of the corresponding cluster
 algebra $\mathcal{A}_{\Omega}$ is given. The images are called preprojective cluster variables.
 It is proved that all preprojective cluster variables other than $u_i$ have denominators
 $u^{\underline{\mbox{dim}}M}$ in their irreducible fractions of integral polynomials,
 where $M$ is the corresponding preprojective module or preinjective module.
 In case the valued quiver $(\G,\Omega)$ is of finite type, the denominator theorem holds with respect to any cluster.
 Namely, let $\underline{x}=(x_1,\cdots, x_n)$ be a cluster of the cluster algebra $\mathcal{A}_{\Omega}$, and $V$ the
 cluster tilting object in $\mathcal{C}(\Omega)$ corresponding to
 $\underline{x}$, whose endomorphism algebra is denoted by $\Lambda$. Then the denominator of any cluster variable $y$ other than $x_i$ is
 $x^{\underline{\mbox{dim}}M}$, where $M$ is the indecomposable $\Lambda-$module corresponding
to $y$.
  This result is a generalization of the corresponding result in [CCS2] to non simply-laced case.

\medskip


\medskip

}
\end{minipage}
\end{center}
\medskip

\begin{center}

\textbf{1. Introduction}\end{center}

\medskip

Clusters and cluster algebras are defined and studied by Fomin and
Zelevinsky [FZ1-2] [BFZ] in order to provide an algebraic framework
for total positivity and canonical bases in semisimple algebraic
groups. Recently there are many works to link representation theory
of quivers with cluster algebras, see amongst others [MRZ], [BMRRT],
[BMR1-2],[BMRT], [CC], [CCS1-2], [K], [CK1-2], [Z1-3]. In [BMRRT], also
[CCS1],  the authors have defined cluster categories and related
cluster tilting theory with clusters and seeds.
 Now cluster categories have become a successful model to understand (acyclic) cluster algebras. For example,
 Caldero and Keller [CK1] [CK2] realize non simply-laced cluster algebras by using cluster categories of quivers via Hall algebra
 approach. From this realization, they proved that, see also [BMRT], Caldero-Chapoton's formula [CC] gives a one-to-one correspondence from
 the set of exceptional objects in
 the cluster category to the set of cluster variables of the corresponding cluster algebra,
 this bijection maps the shift of indecomposable projective modules to the initial cluster variables of the initial seed,
 and sends cluster tilting objects to clusters. Cluster-tilted algebras are by definition, the endomorphism algebras of
 cluster tilting objects in the cluster categories. It provides a class of finite dimensional algebras which is close to the class
 of (quasi-)tilted algebras, but they are different. Their module categories have close connections with cluster categories [BMR1] [Z2] [ABS].

\medskip

The main aims of the paper is the following: The first one is to determine the representation type of cluster-tilted algebras. It is proved that any two cluster-tilted algebras defined in the cluster category $\mathcal{C}(H)$ have same representation type, in particular,
they all have the same representation type as $H$. This result is a consequence of Krause [Kr]
since any two cluster-tilted algebras defined in the same cluster category share a same factor category.
 The second aim is to study the denominators of cluster variables by applying the BGP-reflection functors defined
in [Z1] and their corresponding isomorphisms of cluster algebras in
[Z3]. Analogously with
preprojective or preinjective modules of algebras, we introduce the
notion of preprojective cluster variables. It is proved that the
exponents of the denominators of all preprojective cluster variables
which correspond preprojective modules or preinjective modules $M$
are $\underline{\mbox{dim}}M$. If the cluster algebras is of finite
type, this denominator theorem for cluster variables holds with
respect to any fixed cluster.

\medskip

This paper is organized as follows: In Section 2, some basic notions
 which will be needed later on are recalled. In
Section 3, it is proved that any two cluster-tilted algebras which are defined in the same cluster category share a same factor category, hence they have same representation type.  In Section 4, the notation of
preprojective cluster variables is introduced for any acyclic
cluster algebra. The preprojective cluster variables of a cluster algebra correspond bijectively to the indecomposable preprojective objects of the
corresponding cluster category, i.e. the indecomposable preprojective module,
the indecomposable preinjective objects, or the first shift of indecomposable projective modules in this cluster category. It is proved that the denominator
theorem for all preprojective cluster variables holds with respect to an acyclic seed. For any cluster
algebra of finite type, there is a bijection between the set of indecomposable objects in the corresponding cluster category and the set of cluster variables of this cluster algebras, see Section 4 or [Z1, Z3], this bijection induces a bijection between the set of cluster tilting objects to the set
of clusters. Fix a cluster $\underline{x}=(x_1,
x_2, \cdots, x_n)$, denote by $V$ the (basic) cluster tilting object
corresponding to $\underline{x}$ and by $B$ the corresponding
cluster-tilted algebra. It is proved that any cluster variables $v$ other than
$x_i$ has denominator $x^{\underline{\mbox{dim}}M}$ with respect to the cluster $\underline{x}$, where $M$ is the
 indecomposable $B-$module corresponding to $v$ (compare [FZ2] [CCS1-2]).

 \begin{center}

\textbf{2. Basics on cluster categories and cluster algebras.}
\end{center}

Let $\cal{H}$ be a hereditary abelian category defined over a
field $K$ with finite
dimensional Homomorphism and extension spaces, and with tilting objects. Assume that its
 Grothendieck group is isomorphic to \textbf{$Z^n$} [HRS]. The endomorphism algebra of a tilting object in $\mathcal{H}$ is called a quasi-tilted algebra.
  Denote by $\mathcal{D} =
    D^{b}(\cal{H})$ the bounded derived category of $\cal{H}$ with shift
functor $[1]$. $\mathcal{D}$ has Auslander-Reiten triangles (AR-triangles for short), $\tau$ is the AR-translation. For any category $\cal{T}$, we will denote by
$\ind\cal{T}$ the subcategory of isomorphism classes of
indecomposable objects in $\cal{T}$; depending on the context we
shall also use the same notation to denote the set of isomorphism
classes of indecomposable objects in $\cal{T}$.
\medskip

The orbit category $\mathcal{D} / \tau^{-1}[1]$ is called the
cluster category of type $\mathcal{H}$, which is denoted by
$\mathcal{C(H)}$ ([BMRRT]), see also [CCS1]. It is a triangulated
category with shift functor $[1]$ [Ke] and has Auslander-Reiten triangles, the AR-translation $\tau$ is induced from AR-translation of $\mathcal{D}$ [BMRRT]. When
 $\mathcal{H}$ is the module category of a hereditary algebra $H$, or more general, of a (quasi-)tilted algebra $A$,
or $\mathcal{H}$ is the category of representations of a valued
quiver $(\G, \Omega)$, the corresponding cluster category is denoted
by $\mathcal{C}(H)$, $\mathcal{C}(A)$ or $\mathcal{C}(\Omega)$ respectively. We use $H$ to denote the tensor algebra of the species $\mathcal{M}$
of the valued quiver $(\G, \Omega)$, hence $modH$ is equivalent to the category of representations of the species $\mathcal{M}$ of $(\G, \Omega)$. We denote by $E_1, \cdots, E_n$ complete list of simple objects in $\mathcal{H}$; by $P_1,\cdots,  P_n$ (or $I_1, \cdots , I_n$) the complete list of indecomposable projective (injective respectively) representations in $\mathcal{H}$. Denote $Hom_{\mathcal{C(H)}}(X,Y)$ simply by $Hom(X,Y),$
 and define $Ext^1(X,Y)=Hom(X,Y[1]).$

\medskip

An object $X$ in $\mathcal{C(H)}$ is called exceptional
  if Ext$^1(X,X)=0$. The set of isomorphism
  classes of indecomposable exceptional objects in $\mathcal{C(H)}$ is denote by $\mathcal{E(H)}$.
  If $\mathcal{H}$ is the category of representations of a
  valued quiver $(\G,\Omega)$, this set is denoted by $\mathcal{E}(\Omega)$ sometimes.  An object $V$ in
  $\mathcal{C(H)}$ is called a cluster tilting object if it is exceptional and has $n$ non-isomorphic
  indecomposable direct summands. An exceptional object $M$ in $\cal{C(H)}$ with $n-1$ non-isomorphic direct summands is called almost
  complete.  For a cluster tilting
  object $V$ in $\mathcal{C(H)}$, the
  endomorphism ring End$V$ is called the
  cluster-tilted algebra of $V$ [BMRRT] [BMR1].

\medskip

For a valued graph $\G$, we denote by $\Phi=\Phi^+\bigcup\Phi^-$ the
set of roots of the corresponding Kac-Moody Lie algebra. Let
$\Phi_{\geq -1}=\Phi^+\bigcup\{-\alpha_i\ | \ i=1,\cdots n\ \}$
denote the set of almost positive roots, i.e. the positive roots
together with the negatives of the simple roots.
  Let $s_i$ be the Coxeter generators of the
Weyl group of $\Phi$ corresponding to $i\in \G $. We recall from
[FZ2] that the "truncated reflections" $ \sigma_i$ of $\Phi_{\geq
-1}$ are defined as follows:
$$  \sigma_i(\alpha)=\left\{ \begin{array}{ll} \alpha & \alpha=-\alpha_j,\ j\not=i \\ s_i(\alpha) & \mbox{otherwise.} \end{array}\right.$$

\medskip

 When $\G$ is a Dynkin graph, these truncated reflections were shown in [FZ2] to be one of
the main ingredients of constructions (see also [MRZ]). For Dynkin graph $\G$, the set of vertices can be divided into two
completely disconnected subsets as $\G =\G ^+\sqcup \G ^-$,  and one
can define:
$$\ \ \ \ \sigma_{\pm}=\prod_{i\in \G ^{\pm}}\sigma_i.$$
Then there is a so-called "compatibility degree" $(\ ||\ )$ defined on pairs of almost positive roots of $\G$.
It is uniquely defined by the following two
properties:

 $$\begin{array}{lllcl}& &(-\alpha_i||\beta)&=&\mbox{max}(n_i(\beta),
 0),\\
  &&(\sigma_{\pm}\alpha||\sigma_{\pm}\beta )&=&(\alpha ||\beta),\end{array}$$
  for any $\alpha , \beta \in \Phi_{\geq -1},$ any $i\in \G,$ where $\beta =\sum _in_i(\beta)\alpha_i$ [FZ2].
\medskip

Now we define a map for any valued quiver $(\G, \Omega)$ from
ind$\mathcal{C}(\Omega)$ to $\Phi_{\ge-1}$ as follows: for any $X\in
\mbox{ind}(\mbox{mod}H\vee H[1]),$
$$\gamma_{\Omega}(X)=\left\{
\begin{array}{lrl}\underline{\mbox{dim}}X & \mbox{ if } & X\in \mbox{ind}H;\\
&&\\
-\underline{\mbox{dim}}E_i& \mbox{ if } &X=P_i[1],\end{array}\right.$$ where
$\underline{\mbox{dim}}X$ denotes the dimension vector of the representation
$X$. In general, this map $\gamma_{\Omega}:\
\ind\mathcal{C}(\Omega)\rightarrow
 \Phi_{\geq -1} $ is surjective, but not injective. Denote by $\Phi^{sr}_{\geq-1}$ the set of real Schur roots, i.e. $\underline{dim}X$, where $X$ is exceptional $H-$modules, together with negatives of simple roots. Hence $\gamma_{\Omega}$ induces a bijection from $\mathcal{E}(\Omega)$ to $\Phi^{sr}_{\geq-1}$, which is still
 denoted by $\gamma_{\Omega}$.
\medskip

\medskip

We recall some basic notation on cluster algebras which can be found
in the papers by Fomin and Zelevinsky [FZ1-2].
The cluster algebras we deal with in this paper are defined on a
trivial semigroup of coefficients.
\medskip

 Let
$\mathcal{F}=\mathbf{Q}(u_1,u_2,\cdots, u_n)$ be the field of
rational functions in indeterminates  $u_1,u_2,\cdots, u_n.$
 Set $\underline{u}=(u_1,u_2,\cdots,
u_n).$ Let $B=(b_{ij})$ be an $n\times n$ skew-symmetrizable integer
matrix. A pair $(\underline{x}, B)$, where
$\underline{x}=(x_1,x_2,\cdots, x_n)$ is a transcendence base of
$\mathcal{F}$ and where $B$ is an $n\times n$ skew-symmetrizable
integer matrix, is called a seed. Fix a seed $(\underline{x}, B)$,
$z$ in the base $\underline{x}$. Let $z'$ in
 $\mathcal{F}$ be such that
 $$zz'=\prod_{b_{xz}>0}x^{b_{xz}}+\prod_{b_{xz}<0}x^{-b_{xz}}.$$

 Now, set $\underline{x}':=\underline{x}-\{z\}\bigcup \{z'\}$
  and $B'=(b'_{xy})$ such that

  $$ b'_{xy}=\left\{ \begin{array}{lccccl}-b_{xy}&&&&&\mbox{if } x=z
  \mbox{ or }y=z,\\
  b_{xy}+1/2(|b_{xz}|b_{zy}+b_{xz}|b_{zy}|)&&&&&
  \mbox{otherwise.}\end{array}\right. $$

  The pair $(\underline{x}', B')$ is called the mutation
  of the seed $(\underline{x}, B)$ in direction $z$, it is also a seed.
  The "mutation equivalence $\approx $" is an equivalence relation
  on the set of all seeds generated
  by the mutation.

The cluster algebra $\mathcal{A}_{B}$ associated to the
skew-symmetrizable matrix $B$ is by definition the subalgebra of
$\mathcal{F}$ generated by all $x_i$ in $\underline{x}$ such that
$(\underline{x}, B')\approx(\underline{u}, B).$ Such
$\underline{x}=(x_1,x_2,\cdots, x_n)$ is called a cluster of the
 cluster algebra $\mathcal{A}_{B}$ or simply of $B$, and any
 $x_i$ is called a cluster variable. The set  of all
 cluster variables is  denoted by $\chi _B .$  If the set $\chi _B$ is finite, then the cluster
  algebra $\mathcal{A}_{B}$ is said to be of finite type. Cluster algebras of finite type can be characterized by Dynkin diagrams [FZ2].
\medskip

\medskip

 \begin{center}

\textbf{3. Cluster-tilted algebras.}
\end{center}

  Since any cluster tilting object $V$ in cluster category $\mathcal{C}(\cal{H})$
  is induced by a tilting object in a hereditary abelian category
  $\cal{H}'$, derived
  equivalent to $\cal{H}$, we may assume that,
  without loss the generality (compare [BMRRT] [Z2]),
   $V$
  is a tilting object in $\cal{H}$, and then it is a cluster tilting object in
  $\mathcal{C}(\cal{H})$.
  We have the quasi-tilted algebras $A=\mbox{End}_{\cal{H}}V$ and the
  cluster-tilted algebra $\Lambda=\mbox{End}V$.
\medskip

The Hom functor
 $G=\mbox{Hom}(V,-)$ induces a dense and full functor
 from the
  cluster category $\mathcal{C(H)}$ to $\L-$mod. It induces an equivalence $\overline{G}$ from the factor category
   $ \mathcal{C(H)}/\mbox{add}(\tau
  V)$ to $\L-\mbox{mod}$. This is proved in [BMR1] for the case where $\cal{H}$ is module
category over a hereditary algebra, and generalized to any
hereditary category in [Z2]. Following [CB], [Kr], a ring $A$
  is called generically wild if there is a
  generic $A-$module $M$ such that End$_A(M)$ is not a PI-ring.
\medskip

\medskip

\textbf{Lemma 3.1 [BMR1][Z2].} Let $V$ be a cluster tilting object in
$\cal{C(H)}$ and
 $\Lambda=\mbox{End}V$ the cluster-tilted algebra.
 Then  $\overline{G}: \mathcal{C(H)}/\mbox{add}(\tau
  V)\rightarrow \L-\mbox{mod}$ is an equivalence.

\medskip

From this result, one can compare any two cluster-tilted algebras as following:

\medskip

\textbf{Proposition 3.2.} Let $V$ and $V'$ be cluster tilting objects in $\cal{C(H)}$,
 $\Lambda=\mbox{End}V$ and $\Lambda'=\mbox{End}V'$ the corresponding
 cluster-tilted algebras. Then  $\frac{\Lambda-\mbox{mod}}{\mbox{add}(\mbox{Hom}(V,\tau V'))}
 \approx \frac{\Lambda '-\mbox{mod}}{\mbox{add}(\mbox{Hom}(V',\tau V))}$, and the equivalence
 induces an isomorphism of the AR-quivers between
$\frac{\Lambda-\mbox{mod}}{\mbox{add}(\mbox{Hom}(V,\tau V'))}$ and
$\frac{\Lambda '-\mbox{mod}}{\mbox{add}(\mbox{Hom}(V',\tau V))}$.
Furthermore, $\L$ is generically wild if and only if $\L'$ is
generically wild.
\medskip

\textbf{Proof.} Denote by $G=\mbox{Hom}(V,-).$ The
induced functor $\bar{G}: \mathcal{C(H)}/\mbox{add}(\tau   V)
\rightarrow \L-\mbox{mod}$ is an equivalence by Lemma 3.1. We
consider the composition of functor $\bar{G}$ with the quotient
functor $Q: \L-\mbox{mod} \rightarrow \frac{\Lambda
-\mbox{mod}}{\mbox{add}(\mbox{Hom}(V,\tau V'))}$, which is denoted
by $G_1$. The functor $G_1$ is full and dense since $\bar{G}$ and
$Q$ are. Under the equivalence $\bar{G}$, $\tau V'$
corresponds to Hom$(V, \tau V').$ For any morphism $f: X\rightarrow
Y$ in the category $\frac{\mathcal{C(H)}}{\mbox{add}(\tau V)}$,
$\bar{G}(f): G(X)\rightarrow G(Y)$ factors through
add$(\mbox{Hom}(V,\tau V'))$ if and only if $f$ factors through
add$\tau V'.$ Then $G_1$ induces an equivalence, denoted by
$\bar{G_1},$  from the category $\frac{\cal{C(H)}}{\mbox{add}(\tau
(V\oplus V'))}$ to  the category $\frac{\Lambda
-\mbox{mod}}{\mbox{add}(\mbox{Hom}(V,\tau V'))}.$ It can be proved
in a similar way as Proposition 3.2. in [BMR1] that $\bar{G}$
preserves Auslander-Reiten sequences and then $\bar{G_1},$ induces
an isomorphism of the AR-quivers between
$\frac{\cal{C(H)}}{\mbox{add}(\tau (V\oplus V'))}$ and
$\frac{\Lambda -\mbox{mod}}{\mbox{add}(\mbox{Hom}(V,\tau V'))}.$
Similarly, we have a functor $G'_1$ from
$\frac{\cal{C(H)}}{\mbox{add}(\tau V')}$ to
$\frac{\L'-\mbox{mod}}{\mbox{add}(\mbox{Hom}(V', \tau V))}.$ It
induces an equivalence from the category
$\frac{\cal{C(H)}}{\mbox{add}(\tau (V\oplus V'))}$ to the category $\frac{\Lambda
'-\mbox{mod}}{\mbox{add}(\mbox{Hom}(V',\tau V))}$ and this
equivalence induces an isomorphism of the AR-quivers between
$\frac{\cal{C(H)}}{\mbox{add}(\tau (V\oplus V'))}$ and
$\frac{\Lambda' -\mbox{mod}}{\mbox{add}(\mbox{Hom}(V',\tau V))}.$
Then $\frac{\Lambda-\mbox{mod}}{\mbox{add}(\mbox{Hom}(V,\tau
V'))}\approx \frac{\Lambda '-\mbox{mod}}{\mbox{add(Hom}(V',\tau
V))},$ and the equivalence induces an isomorphism  of the AR-quivers
between $\frac{\Lambda-\mbox{mod}}{\mbox{add}(\mbox{Hom}(V,\tau
V'))}$ and $\frac{\Lambda
'-\mbox{mod}}{\mbox{add}(\mbox{Hom}(V',\tau V))}.$  Since $\L$ and
$\L'$ are artin algebras, it follows from Corollary 3.4. in [Kr]
that $\L$ is generically wild if and only if $\L'$ is generically
wild. The proof is finished.

\medskip

\textbf{Corollary 3.3.}  Let $H$ be a finite dimensional hereditary algebra
over an algebraically closed field, $V$ a cluster tilting object in the
cluster category $\mathcal{C}(H)$  and
$\Lambda=\mbox{End}V$ the corresponding
 cluster-tilted algebra. Then

 (1).  $\L$ is of representation  finite type if and only if $H$ is of representation
finite type; in this case, the numbers of the isomorphism classes of
indecomposable modules over $H$ and over $\L$ respectively are the
same.

(2). $\L$ is tame (or wild) if and only if $H$ is tame (or wild).
\medskip

\textbf{Proof.} Part (1) is proved in Corollary 2.4 in [BMR1]. We prove part (2).
Note that $H$ is a cluster-tilted algebra End$H$,
then by Proposition 3.2., $\L$ is generically wild if and only if
$H$ is generically wild. From Theorem 4.1 in [CB], for a finite dimensional algebra
over an algebraically closed field, it is wild if and only if it is
generically wild. It implies that $\L$ is wild if and only if $H$ is
wild; and $\L$ is tame if and only if $H$ is tame. The proof is
finished.

\medskip

Now we assume that $\bar{V}$ is an almost complete cluster tilting object in
$\mathcal{C(H)}$ and $V_k, \ V_k^*$ are the two complements of
$\bar{V}$, i.e. $V=\bar{V}\oplus V_k$ and  $V'=\bar{V}\oplus V_k^*$
form a tilting mutation at $k$ [BMR2]. Denote the simple top of
projective $\L-$module Hom$(V,V_k)$ by $L_k$ and the simple
top of projective $\L'-$module Hom$(V',V^*_k)$ by $L'_k.$ As a
consequence of Proposition 3.2., we have the following (compare
Theorem 4.2. in [BMR1]):
\medskip

\textbf{Corollary 3.4.} Let $V$ and $V'$ be tilting objects in
$\mathcal{C(H)}$.Suppose $V'$ is obtained from $V$ by a mutation at
$k$, and $\L,$  $\L'$ are the corresponding cluster-tilted algebras.
Then $\frac{\Lambda-\mbox{mod}}{\mbox{add}L_k} \approx \frac{\Lambda
'-\mbox{mod}}{\mbox{add}L'_k}$.
\medskip

\textbf{Proof.} From Proposition 3.2, we have that
$\frac{\Lambda-\mbox{mod}}{\mbox{add(Hom}(V,\tau V^*_k))} \approx
\frac{\Lambda '-\mbox{mod}}{\mbox{add(Hom}(V',\tau V_k))}.$ It is
easy to see Hom$(V,\tau V^*_k)\cong L_k$ and Hom$(V',\tau V_k)\cong
L'_k$ (compare Lemma 4.1 in [BMR1]), which implies the equivalence
what we want. The proof is finished.

\medskip

 In the rest of this section, we assume the hereditary category $\mathcal{H}$ is
the category of finite dimensional left module of a finite
dimensional algebra $H$ over a filed $K$. Assume that $V$ is a tilting $H-$module. Then Hom$_H(V,-)$ induces a derived
equivalence from $D^b(H)$ to $D^b(A)$, where $A$ is the tilted
algebra End$_HV$. It induces a triangle equivalence from the cluster
category $\mathcal{C}(H)$ to the cluster category $\mathcal{C}(A)=D^b(A)/\tau^{-1}[1]$.
This equivalence is denoted by $R_V$ [Z1, Z2].
\medskip

\textbf{Proposition 3.5.} Let $H$, $V$, $A$ be as above and $\L=\mbox{End}V$ the
cluster-tilted algebra. Then we have the commutative diagram: \[ \begin{CD}
\mathcal{C}(H) @>R_V
>> \mathcal{C}(A)\\
@V VV  @VVV  \\
\frac{H-\mbox{mod}}{\mbox{add}(\tau V)} @>\approx>> \frac{\L-\mbox{mod}}{\mbox{add}(\mbox{Hom}(V, \tau H))},
\end{CD} \]  where the left vertical map is $Q_H\circ \mbox{Hom}(H,-),$  the right vertical map is $Q_{\L}\circ \mbox{Hom}(A,-)$, and $Q_H: H-\mbox{mod}\rightarrow \frac{H-\mbox{mod}}{\mbox{add}\tau V}$,  $Q_{\L}: \L-\mbox{mod}\rightarrow \frac{\L-\mbox{mod}}{\mbox{add}(\mbox{Hom}(V, \tau H))}$ are the quotient maps.

\medskip

\textbf{Proof.} $R_V(V)=A$ is a cluster tilting object in $\mathcal{C}(A)$, and Hom$(V,-)=\mbox{Hom}(A,-)\circ R_V.$ Then the commutativity of the diagram above follows from Proposition 3.2.

\medskip

For any tilting $H-$module $V$, from Lemma 3.1., we have that
ind$\mathcal{C}(H)=\mbox{ind}\L \vee \{V_i[1]\ | \ i=1, \cdots n\}.$
We  define a map from ind$\mathcal{C}(H)$ to $\mathbf{Z^n}$ by
generalizing the map $\gamma_{H}.$

\medskip

\textbf{Definition 3.6.} Let  $\gamma_{V}:
\mbox{ind}\mathcal{C}(H)\rightarrow \mathbf{Z^n}$ be defined as
follows: for any $X\in \mbox{ind}\L \vee \{ V_i[1] | i=1,\cdots n
\},$

$$\gamma_{V}(X)=\left\{
\begin{array}{lrl}\underline{\mbox{dim}}X & \mbox{ if } & X\in \mbox{ind}\L ;\\
&&\\
-\underline{\mbox{dim}}E_i& \mbox{ if } &X=V_i[1],\end{array}\right.$$

 $\gamma_{V}$ is called
the dimension vector map of the cluster category $\mathcal{C}(H)$
associated to the tilting object $V.$  If $X\in \mbox{ind}\L$, then we
say that $\gamma_{V}(X)>0.$  Denote  by  $\Phi_V=\{ \gamma_{V}(X)\ |
\ X\in \mbox{ind}\mathcal{C}(H)\ \}$ and $\Phi_V^{sr}=\{
\gamma_{V}(X)\ | \ X\in \mathcal{E}(H)\ \}.$

\medskip

In general, this map $\gamma_{V}:\
\ind\mathcal{C}(H)\rightarrow
 \Phi_V $ is not injective. When $V=H$, $\gamma_{V}=\gamma_{H}$ which maps from ind$\mathcal{C}(H)$ to $\Phi_{\ge -1}.$

\medskip

\textbf{Lemma 3.7.} Suppose $V=\oplus _{i=1}^{i=n}V_i$ is a cluster tilting object, $M\in \mbox{ind}\mathcal{C}(H)$ and
$\tau M \rightarrow \oplus_j N_j\rightarrow M\rightarrow \tau M[1]$
is the AR-triangle ending at $M$, where $N_j$ is indecomposable in
$\mathcal{C}(H)$. Then

 (1). if $M\cong V_i$ for some $i$, we
have $\gamma_{V}(M)=\underline{\mbox{dim}}E_i+\sum_{\gamma
_V(N_j)>0} \gamma _V(N_j)$;

(2). if $M\cong V_i[1]$ for some $i$, we have $\gamma_{V}(\tau
M)=\underline{\mbox{dim}}E_i+\sum_{\gamma _V(N_j)>0} \gamma
_V(N_j)$;

 (3). if $M\ncong V_i$ and $M\ncong V_i[1]$ for any $i$, we have $\gamma_{V}(M)+\gamma_{V}(\tau M)=\sum_{\gamma _V(N_j)>0} \gamma _V(N_j)$.

\medskip

\textbf{Proof.} By [BMR1], the sink map and the source map for any
indecomposable $\L-$module $M$ in $\L-$mod are induced from the
corresponding maps in the cluster category $\mathcal{C}(H).$
Then the AR-triangle  $\tau M \rightarrow \oplus N_j\rightarrow
M\rightarrow \tau M[1]$ in $\mathcal{C}(H)$ induces \begin{enumerate} \item the AR-sequence:
 $0\rightarrow \mbox{Hom}(V,\tau M) \rightarrow \oplus_{\scriptsize{\mbox{Hom}(V, N_j)\not=0}} \mbox{Hom}(V, N_j)
 \rightarrow  \mbox{Hom}(V,M)\rightarrow 0$ provided $\mbox{Hom}(V,M)\not=0$  and  is not any projective $\L-$module; or
 \item the source map: $ \mbox{Hom}(V,\tau M) \rightarrow
\oplus_{\scriptsize{\mbox{Hom}(V, N_i)\not=0}} \mbox{Hom}(V, N_i)$ provided  $\mbox{Hom}(V,\tau M)$ is an injective $\L-$module or
\item a
sink map:  $\oplus_{\scriptsize{\mbox{Hom}(V, N_j)\not=0}} \mbox{Hom}(V,
N_j)\rightarrow  \mbox{Hom}(V,M)$ provided $\mbox{Hom}(V, M)$ is
a projective $\L-$module.\end{enumerate} Then by the definition of $\gamma_V$, we have
that (3). $\gamma_{V}(M)+\gamma_{V}(\tau M)=\sum_{\gamma _V(N_j)>0}
\gamma _V(N_j)$ if  $M\ncong V_i$ and $M\ncong V_i[1]$ for any $i$;
 (2).$\gamma_{V}(\tau M)=\underline{\mbox{dim}}E_i+\sum_{\gamma
_V(N_j)>0} \gamma _V(N_j)$ if $M\cong V_i[1]$ for some $i$; and (3).
$\gamma_{V}(M)=\underline{\mbox{dim}}E_i+\sum_{\gamma _V(N_j)>0}
\gamma _V(N_j)$ if $M\cong V_i$ for some $i$.  The proof is
finished.
\medskip

\textbf{Definition 3.8.} Let $V$ be a cluster tilting object in
$\mathcal{C}(\Omega).$ We Define $\sigma_{V}:
\Phi_{\ge-1}^{sr}\rightarrow \Phi_V ^{sr}$  as the map
$\gamma_{\Omega}(X)\mapsto \gamma_{V}(X)$ for any $X\in
\mathcal{E}(H).$
\medskip

Since the set of indecomposable exceptional objects in
$\mathcal{C}(H)$ consists of indecomposable exceptional
$H-$modules and the shifts by $1$ of indecomposable projective
$H-$modules, they are determined by their dimensional vectors. Then
the map $\sigma_V$ is well-defined. When $V$ is a BGP-tilting or
APR-tilting module (at a sink or a source $k$),
$\sigma_{V}=\sigma_{k}.$

\medskip

Now we give an interpretation of Proposition 3.5. in terms of
dimensional vector maps.

\medskip

\textbf{Proposition 3.9.} Let $V$ be a tilting representation of $(\G,\Omega)$, whose endomorphism algebra is denoted by $A$. Then we have the commutative diagram: \[ \begin{CD}
\mathcal{E}(\Omega) @>R_V
>> \mathcal{E}(A)\\
@V\gamma_{\Omega} VV  @VV\gamma_{V}V  \\
\Phi_{\ge-1}^{sr} @>\sigma_{V}>> \Phi_V ^{sr}
\end{CD} \].
\textbf{Proof.} Since $R_V$ is a triangle equivalence from $\mathcal{C}(\Omega)$ to the cluster category $\mathcal{C}(A)$ of $A$, it maps exceptional objects to exceptional objects. Hence $R_V$ induces bijection from $\mathcal{E}(\Omega)$ to $\mathcal{E}(A)$. $A=R_V(V)$ is a cluster tilting object and $\L=\mbox{End}V\cong \mbox{End}A$.  By definition of dimensional vectors map, $\gamma_{A}=\gamma_{V}.$ Then for any indecomposable exceptional object $X\in \mathcal{C}(\Omega),$  $\gamma_{V}R_V(X)=\gamma_{V}(X)=\sigma_{V}(\gamma_{\Omega}(X)).$ the proof is finished.

\medskip

\begin{center}

\textbf{4. Preprojective cluster variables}\end{center}

\medskip

In this section, we assume that $\mathcal{H}$ is the category of
finite dimensional representations over a field $K$ of a species $\mathcal{M}$ of the valued quiver $(\G, \Omega)$.
 We always assume that $\G$ contains no vertex loops and the quiver $(\G, \Omega)$
  contains no oriented cycles. Let $H$ denote the tensor algebra of the species $\mathcal{M}$. It is a hereditary finite dimensional $K-$algebra.

   For any  vertex $k\in \G, $ $s_k\Omega$ denotes the new orientation of $\G$
   by reversing the direction of arrows
  along all edges containing $k$ in $(\G, \Omega)$. A vertex $k\in \G$ is called a sink (or a source) with respect to $\Omega$
  if there are no arrows
  starting (or ending) at vertex $k$.

  \medskip

We recall that $P_i$ (or $I_i)$ is the projective (injective resp.)
indecomposable $H-$module  corresponding to the vertex $i\in \G$,
and $E_i$ are the corresponding simple
 module. If $k$ is a sink (or a source), then
$P_k=E_k$ (resp. $I_k=E_k$) is a simple projective (resp. injective)
$H-$module.
 \medskip

  Let $V=\oplus _{i\in \G-\{k\}}P_i \oplus \tau ^{-1}P_k$. If $k$ is a sink, then $V$ is a tilting $H-$module which is called a BGP-
or an APR-tilting module, $ \mbox{Hom}_H(V,-)$ is the titling functor induced by $V$, which is called a BGP-reflection functor and
is denoted by $S^+_k.$ The following theorem was proved in [Z1] [Z2]
(in a more general case).

\medskip

 \textbf{Theorem 4.1.}[Z1, Z2] For any  sink (or a source) $k$
  of a valued quiver $(\G, \Omega),$ the BGP-reflection functor $S^+_k$ (resp. $S^-_k$) induces
  a triangle equivalence $R(S^+_k)$
  (resp.,$R(S^-_k)$) from
  $\mathcal{C}(\Omega)$ to
  $\mathcal{C}(s_k\Omega)$ and $\gamma_{s_k\Omega}(R(S^+_k)(X))=\sigma_k(\gamma _{\Omega}(X)).$
   Moreover  $R(S^+_k)$ induces a bijection from $\mathcal{E}(\Omega)$ to $\mathcal{E}(s_k\Omega).$

\medskip

\textbf{Definition 4.2.} The triangle equivalence $R(S^+_k)$
from $\mathcal{C}(\Omega)$ to
  $\mathcal{C}(s_k\Omega)$ induced from the reflection functor $S^+_k$
  is called the BGP-reflection functor in $\mathcal{C}(\Omega)$
  at the sink $k$, which is denoted simply by $R^+_k.$ Dually
  for a source $k$, we have the reflection functor $R^-_k$ from
  $\mathcal{C}(\Omega)$ to
  $\mathcal{C}(s_k\Omega)$.
\medskip

\medskip

  Let $k_1, \cdots, k_n$ be an admissible sequence of sinks of $(\G,\Omega)$. For simplicity, we assume that $1, 2, \cdots n$ is
  such an admissible sequence of sinks of $(\G,\Omega).$ Set $C^+=R^+_{n}\cdots R^+_{2}R^+_{1},$
the composition of $R^+_{i}$. $C^+$ is a self-equivalence of
$\mathcal{C}(\Omega)$, it is called the Coxeter functor in the
cluster category $\mathcal{C}(\Omega)$ in [Z3]. For simplicity, we denote
$C^+$ by $C$. The inverse $C^-$ of $C$, which is also called the Coxeter
functor in $\mathcal{C}(\Omega)$ is $C^-=R^-_{1}R^-_{2}\cdots
R^-_{n}.$ For any indecomposable object $X$ in $\mathcal{C}(\Omega)$, $C(X)\cong \tau X$, where $\tau$ is the Auslander-Reiten translation in $\mathcal{C}(\Omega)$. Let $\mathcal{P}$(or $\mathcal{I})$ denote the
set of isomorphism classes of  indecomposable preprojective (preinjective resp.) $H-$modules, $\mathcal{R}$ the set of isomorphism classes of indecomposable regular $H-$modules.
 We denote the union $\mathcal{P}\vee\{ P_i[1]\ |\ i=1, \cdots n \}
   \vee\mathcal{I}$ by $\mathcal{PI}(\Omega).$  Note that
   $$\mathcal{PI}(\Omega)=\{ C^{m}P_k[1]\ \ |
   \  m\in \mathbf{Z}, 1\le k\le n\ \}.$$
The objects in $\mathcal{PI}(\Omega)$ are called preprojective objects in $\mathcal{C}(\Omega).$
\medskip

   If $\G$ is a Dynkin diagram, then ind$\mathcal{C}(\Omega)
   =\mathcal{PI}(\Omega),$ otherwise ind$\mathcal{C}(\Omega)
   =\mathcal{PI}(\Omega)\vee \mathcal{R}.$

   \medskip

\medskip

\textbf{Definition 4.3.} Let $A=(a_{ij})$ be a generalized Cartan matrix corresponding to the diagram $\G$.  For any
index $i$, we define an automorphism $T_i$ of
$\mathcal{F}=\mathbf{Q}(u_1,\cdots,u_n)$ by defining the images of
the indeterminates $u_1, \cdots, u_n$ as follows:

$$T_i(u_j)=\left\{ \begin{array}{lccccl}u_j&&&&&\mbox{if } j\not=i,\\
&&&&&\\
  \frac{\prod_{a_{i k}<0}u_k^{-a_{i k}}+1}{u_i}&&&&&
 \mbox{if } j=i.\end{array}\right.$$

 It is easy to check that all $T_i$ are involutions of $\mathcal{F}$, i.e.
  $T_i^2=\mbox{id}_{\mathcal{F}}.$

Let $B$ be a skew-symmetrizable  matrix corresponding to the valued
quiver $(\G,\Omega)$. Then the Cartan counterpart of $B$ is $A$
[FZ2].

\medskip

 From the definition, all $T_i$ are independent of orientations of
the valued quiver $(\G,\Omega)$ corresponding to $B$, and depends
only on the valued graph $\G$. Let $k$ be a sink in $(\G,\Omega )$.
By Theorem 4.7 in [Z3], $T_k$ sends cluster variables and
 clusters in $\chi_{\Omega }$ to those in
$\chi_{s_k\Omega}$ respectively and  $T_k$ induces an isomorphism
from the cluster algebra $\mathcal{A}_{\Omega}$ to
$\mathcal{A}_{s_k\Omega}$ ($T_k$ induces a so-called strongly
isomorphism from $\mathcal{A}_{\Omega}$ to
 $\mathcal{A}_{s_k\Omega}$).
\medskip

 Under the assumption that the sequence $1,2,\cdots, n$ is an admissible sequence of sinks of $(\G,\Omega)$, we define an automorphism $T_{\Omega}$
  of $\mathcal{F}$ as $T_{\Omega}=T_{n}\cdots T_{2}T_{1}.$ $T_{\Omega}$ induces an automorphism of cluster algebra
  $\mathcal{A}_{\Omega}$.   $T_{\Omega}$ and its inverse
$T^{-}_{\Omega}$ are called the Coxeter automorphisms of the cluster
algebra  $\mathcal{A}_{\Omega}$. $T_{\Omega}$ is simply denoted by
$T$.
\medskip

Set:
  $$\chi'_{\Omega}=\{ T^m(u_k)\  |  \  m\in \mathbf{Z}, 1\le k \le
  n \  \}.$$
\medskip

  When $\G$ is of finite type, we have that $\chi'_{\Omega}=\chi_{\Omega}$
  is the set of cluster variables of $(\G, \Omega)$ and $C^m(P_i[1])\mapsto T^m(u_i)$ is a bijection from $\mathcal{C}(\Omega)$ to $\chi_{\Omega}$
  which sends $P_i[1]$ to $u_i$ and sends tilting objects to clusters by Corollary 3.4.
  and Theorem 4.7. in [Z3].  For any valued quiver, it is proved in Theorem 4.13 of [Z3] that  all elements $T^m(u_i)$
  in $\chi'_{\Omega}$ are cluster variables of cluster algebra $\mathcal{A}_{\Omega}$.
   The main aim in this section is to prove that the denominators of these clusters
$T^m(u_i)$ as the reduced fractions of integral polynomials are
$u^{\gamma_{\Omega}(C^m(P_i[1]))}$. As an immediately consequence, the map
 $\phi_{\Omega}: C^k(P_i)\mapsto T^k(u_i) ,$ for any $k\in\mathbf{Z}, \ i=1,\cdots, n, $ is
a bijection from $\mathcal{PI}(\Omega)$ to $\chi'_{\Omega}$.
Clusters and cluster variables in $\chi'_{\Omega}$ will be called
preprojective clusters and preprojective cluster variables
respectively.

  \medskip

Before we state the theorem, we recall a useful notation from [BMRT]. From [FZ1], any cluster variable $x$ ($\not= u_j \forall j$) is a Laurent polynomial, i.e. it can be written as
a reduced fraction $f/m$ where $f$ is an integer polynomial in $u_1, \cdots, u_n$ and $m=u^{\underline{a}}=u_1^{a_1}u_2^{a_2}\cdots u_n^{a_n}$, where $\underline{a}=(a_1,\cdots , a_n)$ is a non-negative vector, is a monomial in $u_1, \cdots, u_n$. Following a technique used in [BMRT], we call a polynomial $f$ "positive" provided $f(e_i)>0$ for any $e_i=(1,\cdots, 1,0,1,\cdots, 1)$ for $i=1,\cdots n.$  It is clear that if $x=f/m$ with a "positive" polynomial $f$ and a monomial $m$, then $x$ is in reduced form [BMRT]. In this case or $x=u_i=\frac{1}{u_i^{-1}}$ for some $i$,  $x$ is called of the good reduced form.
\medskip

\textbf{Theorem 4.4.} Any element $T^m(u_k)$ in $\chi'_{\Omega}$ is
a cluster variable of the cluster algebra $\mathcal{A}_{\Omega}$
with denominator $u^{\gamma_{\Omega}(C^m(P_k[1]))}$. Furthermore,
the assignment $\phi_{\Omega}: C^m(P_k[1])\mapsto T^m(u_k), \forall
m\in \mathbf{Z}, \  k\in \G ,$ is a bijection from
$\mathcal{PI}(\Omega)$ to $\chi'_{\Omega}$ such that  $P_k[1]$
corresponds to $u_k$.
\medskip

\textbf{Proof.} Let $B$ be the skew-symmetrizable  matrix
corresponding to the valued quiver $(\G, \Omega)$. If $k$ is a sink
or a source, we denote by $s_kB$ the skew-symmetrizable  matrix
corresponding to the valued quiver $(\G, s_k\Omega).$  As above we
assume that $1, 2, \cdots, n$ is an admissible sequence of sinks of
$(\G, \Omega)$. By definition, $T_1$ sends the seed $((u_1,\cdots ,
u_n), B)$ (which is denoted by $(\underline{u}, B)$ for simplicity)
of cluster algebra $\mathcal{A}_{\Omega}$ to the seed
$(T_1(\underline{u}),B)=((T_1(u_1),u_2, \cdots , u_n), B)$ of
the cluster algebra $\mathcal{A}_{s_k\Omega},$ which can be viewed as
one obtained by seed mutation in direction $1$ from the seed $(\underline{u}, s_1B).$  Similarly, $(T_2T_1(\underline{u}), B)$ is a
seed of $\mathcal{A}_{s_2s_1\Omega}$. By induction,
$(T(\underline{u}), B)$ is a seed of $\mathcal{A}_{\Omega}$ since  $s_n\cdots s_1\Omega=\Omega.$ In this
way, we prove that for any non-negative integer $m$,
$(T^m(\underline{u}), B)$ is a seed of $\mathcal{A}_{\Omega}.$ For $m<0$, we use $T^{-1}=T_1\cdots T_n$ to
replace $T$, the same argument implies $(T^{-m}(\underline{u}), B)$
is a seed of $\mathcal{A}_{\Omega}$ for any positive
integer $m$. This proves the first statement in the theorem.

   In the following, we will prove that all $T^m(u_i)$  can be written as a
good reduced form $\frac{f}{u^{\gamma_{\Omega}(C^m(P_i[1]))}}$ with
$f$ a "positive" polynomial in $u_1,\cdots, u_n$. Firstly we note
that $u_i=\frac{1}{u_i^{-1}}=\frac{1}{u^{\gamma_{\Omega}(P_i[1])}}$
is a good reduced form;  and if $\G$ is of finite type, then
$\phi_{\Omega}$ is a bijection which sends tilting objects to
clusters, and sends $P_i[1]$ to $u_i$ (compare Theorem 4.7 and Remark 4.9 in [Z3]). Secondly, we
note that if $P_i[1]$ is a direct summand of a tilting object, i.e.
$V=\bar{V}\oplus P_i[1]$, then Hom$(P_i,\bar{V})=0$, this means the
$i-$th component of vector $\gamma_{\Omega}(\bar{V})$ is zero. Now
we start at the $n-$tuple $(u_1,\cdots,u_{n-1}, T(u_n))$.
$((u_1,\cdots,u_{n-1}, T(u_n)), s_nB)$ is a seed obtained from the
initial seed  $(\underline{u}, B)$ by seed mutation in the direction
$n$, where $n$ is a source in $(\G, \Omega).$ We note that
$T(u_n)=T_n(u_n)$.
 Then $T(u_n)u_n=\prod_{k<n}u_k^{-a_{n,k}}+1$, i.e. $T(u_n)=\frac{\prod_{k<n}u_k^{-a_{n,k}}+1}{u_n}$
which is a good reduced form
$\frac{f_1}{u^{\gamma_{\Omega}(C^1(P_n[1]))}}$ with that $f_1$ is a
"positive" polynomial in $u_1,\cdots, u_n$ since $C^1(P_n[1])$ is a
simple injective module. Similarly, $((u_1,\cdots,u_{n-2},
T(u_{n-1}),  T(u_n)), s_{n-1}s_nB)$ is a seed obtained from the
the seed $((u_1,\cdots,u_{n-1}, T(u_n)), s_nB)$ by seed mutation in
the direction $n-1$, where $n-1$ is a source in $(\G, s_{n}\Omega).$ The reason is the following: Firstly,
 $T_{n}T_{n-1}(u_{n-1})=T(u_{n-1})$ from the definition of $T_i$; Secondly, from the definition of $T_{n-1}$ we have $T_{n-1}(u_{n-1})u_{n-1}=\prod_{k\not= n-1}u_k^{-a_{n-1,k}}+1$
and then
$T(u_{n-1})=T_nT_{n-1}(u_{n-1})=\frac{T_n(u_n)^{-a_{n-1,n}}\prod_{k<n-1
}u_k^{-a_{n-1,k}}+1}{u_{n-1}}$. Since $C^1(P_{n-1}[1])$ is an
injective $H-$module and the AR-triangle starting at
$C^1(P_{n-1}[1])$ is $$C^1(P_{n-1}[1])\rightarrow
C^1(P_{n}[1])^{-a_{n-1,n}}\oplus \oplus_{i\not=n,
n-1}P_i[1]^{-a_{n-1,i}}\rightarrow P_{n-1}[1]\rightarrow
C^1(P_{n-1}[1])[1],$$ by Lemma 3.8, we have that
$\gamma_{\Omega}(C^1(P_{n-1}[1]))=
\underline{\mbox{dim}}E_{n-1}+(-a_{n-1,n})\gamma
_{\Omega}(C^1(P_n[1])).$ Then
$T(u_{n-1})=\frac{f_2}{u^{\gamma_{\Omega}(C^1(P_{n-1}[1]))}}$ where
$f_2=f_1^{-a_{n-1,n}}\prod_{k<n-1}u^{-a_{n-1,k}}+u^{\gamma_{\Omega}(C^1(P_n[1]))}$ is a
"positive" polynomial since $f_1$ is "positive".
  By induction on $k$, we get $T(u_k) =\frac{f_k}{u^{\gamma_{\Omega}(C^1(P_k[1]))}}$
with a "positive" polynomial $f_k$ for all $k$. So we have that
$(T(\underline{u}), B)$ is a seed with that all $T(u_i)$ can be
written as a good reduced form with denominator
$u^{\gamma_{\Omega}(C^1(P_i[1]))}$.

Now we assume that  $((T^{k}(u_1),\cdots, T^{k}(u_i),
T^{k+1}(u_{i+1}),  \cdots, T^{k+1}(u_n)),\\  s_{i+1}\cdots s_nB),$
where $k\ge1, $  is a seed with that all cluster variables
$T^{k}(u_j), \\  \forall j\le i$  and  $T^{k+1}(u_{j})\ \forall j>i$
 are of good reduced forms with denominators
$u^{\gamma_{\Omega}(C^k(P_j[1]))}$ or
$u^{\gamma_{\Omega}(C^{k+1}(P_j[1]))}$ respectively. We will show
that the cluster variable $T^{k+1}(u_i)$ can be written as a good
reduced form with denominator
$u^{\gamma_{\Omega}(C^{k+1}(P_i[1]))}.$

  As before, the $n+1-$tuple $((T^{k}(u_1),\cdots,T^k(u_{i-1}), T^{k+1}(u_i), T^{k+1}(u_{i+1}),\\
    \cdots, T^{k+1}k(u_n)), s_is_{i+1}\cdots s_nB)$ is a pair obtained from the seed $((T^{k}(u_1),\cdots, T^{k}(u_i),\\ T^{k+1}(u_{i+1}),
      \cdots, T^{k+1}(u_n)), s_{i+1}\cdots s_nB)$ by seed mutation in the direction $i$, hence it is a seed.
The reason is the following: From the definition of $T_i$, we have
$T_i(u_i)u_i= \prod_{j\not=i}u_k^{-a_{ik}}+1$ which can be written
as $$T_iT_{i-1}\cdots T_1(u_i)u_i= \prod_{j>i}T_iT_{i-1}\cdots
T_1(u_j)^{-a_{ij}}\prod_{j<i}u_j^{-a_{i j}}+1.$$  By applying
$T_n\cdots T_{i+1}$ to the two sides of the equality, we get
$$T(u_i)u_i= \prod_{j>i}T(u_j)^{-a_{i j}}\prod_{j<i}u_j^{-a_{i j}}+1.$$ and
then  $$T^{m+1}(u_i)T^m(u_i)= \prod_{j>i}T^{m+1}(u_j)^{-a_{i j}}\prod_{j<i}T^m(u_j)^{-a_{i j}}+1.$$

      We have the following AR-sequence in $H-$mod: $$\begin{array}{l}0\rightarrow C^{m+1}(P_i[1])\rightarrow \oplus_{j>i}C^{m+1}(P_j[1])^{-a_{ij}}
\oplus \oplus _{j<i}C^m(P_j[1])^{-a_{ij}}\rightarrow C^m(P_i[1])\rightarrow 0.\end{array}$$
Then $\gamma_{\Omega}(C^{m+1}(P_i[1]))=\sum _{j>i} (-a_{ij})\gamma_{\Omega}(C^{m+1}(P_j[1]))+\sum _{j<i} (-a_{ij})
\gamma_{\Omega}(C^{m}(P_j[1]))-\gamma_{\Omega}(C^{m}(P_i[1])).$
   By assumption that $T^{m+1}(u_j)=\frac{f_j}{u^{\gamma_{\Omega}(C^{k+1}(P_j[1]))}}\  \forall j>i$
and $T^{m}(u_j)=\frac{f_j}{u^{\gamma_{\Omega}(C^{k}(P_j[1]))}} \
\forall j\le i$, all forms are good reduced forms with "positive"
polynomial $f_j$ respectively. Then we have that
$$\begin{array}{rl}T^{m+1}(u_i)=& \frac{\prod_{j>i}(f_j/u^{\gamma_{\Omega}(C^{k+1}(P_j[1]))})^{-a_{i j}}
\bullet \prod_{j<i}(f_j/u^{\gamma_{\Omega}(C^{k}(P_j[1]))})^{-a_{i j}}+1}{f_i/u^{\gamma_{\Omega}(C^{k}(P_i[1]))}}\\
\\
=& \frac{\prod_{j\not=i}(f_j)^{-a_{i j}}+\prod_{j< i}u^{(-a_{i j})\gamma_{\Omega}(C^{k}(P_j[1]))}\prod_{j> i}u^{(-a_{i j})
\gamma_{\Omega}(C^{k+1}(P_j[1]))}}{f_i}\\
&\times \frac{u^{\gamma_{\Omega}(C^{k}(P_i[1]))}}{\prod_{j<
i}u^{(-a_{i j})\gamma_{\Omega}(C^{k}(P_j[1]))} \prod_{j> i}u^{(-a_{i j})\gamma_{\Omega}(C^{k+1}(P_j[1]))}}. \end{array}$$ The first
factor is  a polynomial by the "Laurent Phenomenon" [ZF1], which is
easily to see to be "positive" since all $f_i$ are "positive";  and
the second factor is a monomial $$u^{\sum _{j>i}
(-a_{i j})\gamma_{\Omega}(C^{m+1}(P_j[1]))+
\sum_{j<i}(-a_{i j})\gamma_{\Omega}(C^{m}(P_j[1]))-\gamma_{\Omega}(C^{m}(P_i[1]))}$$
which is $u^{\gamma_{\Omega}(C^{m+1}(P_i[1]))}.$ This proves that
$T^{m+1}(u_i)$ is of a good reduced form with denominator
$u^{\gamma_{\Omega}(C^{m+1}(P_i[1]))}.$

  If $\G$ is a Dynkin diagram, there are two additional cases which may occur:
\\

  \noindent (1). $T^{m+1}(u_i)=u_j$ for some $j$. In this case, $ T^{m+1}(u_i)$ is a good
reduced form $\frac{1}{u^{\gamma_{\Omega}(P_j[1])}}$ since
$\phi_{\Omega}$ is a bijection sending $P_j[1]$ to $u_j$.
\\

\noindent (2). Some $u_j$ appears in the $n-$tuple
$(T^{k}(u_1),\cdots, T^{k}(u_i), T^{k+1}(u_{i+1}),
      \cdots, \\ T^{k+1}(u_n)) $ as components. For simplicity, we assume that $T^k(u_{j_1})=u_{t_{j_1}}$
      for $j_1\in J_1\subset \{ 1, \cdots i\}$ and
      $T^{k+1}(u_{j_2})=u_{t_{j_2}}$ for $j_2\in J_1\subset \{ i+1, \cdots n\}$. The AR-triangle ending at $C^m(P_i[1])$ is
      $$\begin{array}{l} C^{m+1}(P_i[1])\rightarrow \oplus_{j>i}C^{m+1}(P_j[1])^{-a_{ij}}
\oplus \oplus _{j<i}C^m(P_j[1])^{-a_{i j}}\rightarrow C^m(P_i[1])\\
\rightarrow C^{m+1}(P_i[1])[1].\end{array}$$  This triangle induces the AR-sequence  in $H-$mod (compare Lemma 3.7):
 $$\begin{array}{l}0\rightarrow C^{m+1}(P_i[1])\rightarrow \oplus_{j>i, j\notin J_2}C^{m+1}(P_j[1])^{-a_{i j}}
\oplus \oplus _{j<i, j\notin J_1}C^m(P_j[1])^{-a_{i j}}\\
\rightarrow C^m(P_i[1])
\rightarrow 0.\end{array}$$  Then $$\begin{array}{ll}\gamma_{\Omega}(C^{m+1}(P_i[1]))&=\sum _{j>i, j
\notin J_2} (-a_{i j})\gamma_{\Omega}(C^{m+1}(P_j[1]))\\
&+\sum _{j<i, j\notin J_1} (-a_{i j})
\gamma_{\Omega}(C^{m}(P_j[1]))-\gamma_{\Omega}(C^{m}(P_i[1])).\end{array}$$
As before we have the formula for $T^{m+1}(u_i)$:
\small{$$\begin{array}{l}T^{m+1}(u_i)= \frac{\prod_{j>i, j\notin J_2}
(f_j/u^{\gamma_{\Omega}(C^{k+1}(P_j[1]))})^{-a_{i j}}\bullet \prod_{j<i, j\notin J_1}(f_j/u^{\gamma_{\Omega}(C^{k}(P_j[1]))})^{-a_{i j}}
\prod_{j\in J_1 \cup J_2}u_{t_j}^{-a_{i j}}+1}{f_i/u^{\gamma_{\Omega}(C^{k}(P_i[1]))}}\\
\\
= \frac{\prod_{j\notin J_1 \cup J_2, j\not=i}(f_j)^{-a_{i j}}\prod_{j\in J_1 \cup J_2}u_{t_j}^{-a_{i j}}   +\prod_{j< i, j
\notin J_1}u^{(-a_{i j})\gamma_{\Omega}(C^{k}(P_j[1]))}\prod_{j> i, j\notin J_2}u^{(-a_{i j})\gamma_{\Omega}(C^{k+1}(P_j[1]))}}{f_i}\\
\times \frac{u^{\gamma_{\Omega}(C^{k}(P_i[1]))}}{\prod_{j< i,
j\notin J_1}u^{(-a_{i j}) \gamma_{\Omega}(C^{k}(P_j[1]))}\prod_{j>
i, j\notin J_2}u^{(-a_{i j})\gamma_{\Omega}(C^{k+1}(P_j[1]))}}.
\end{array}$$}

The second factor on the right hand of the last equation is the
monomial $u^{\gamma_{\Omega}(C^{m+1}(P_i[1]))}$. The first factor is
a polynomial by the "Laurent Phenomenon" [ZF1], which is "positive"
since all $f_i$ are "positive"  and the monomials $\prod_{j\in J_1
\cup J_2}u_{t_j}^{-a_{i j}}$ and $\prod_{j< i, j\notin J_1}u^{(-a_{i j})\gamma_{\Omega}(C^{k}(P_j[1]))}\prod_{j> i, j\notin J_2}u^{(-a_{i j
})\gamma_{\Omega}(C^{k+1}(P_j[1]))}$ are co-prime.  The last
statement follows from    that $V=\oplus _{j\le i}C^m(P_j[1])\oplus
\oplus_{j>i}C^{m+1}(P_j[1])$ is a cluster tilting object in
$\mathcal{C}(\Omega)$, for $j\in J_1\cup J_2$, $P_{t_j}[1]$ is a
direct summand of $V$, then the $j-$th components of
$\gamma_{\Omega}(C^{m+1}(P_j[1]))$ and of
$\gamma_{\Omega}(C^{m}(P_j[1]))$ are zero. This proves that
$T^{m+1}(u_i)$ is a good reduced form with denominators
$u^{\gamma_{\Omega}(C^{m+1}(P_i[1]))}$.

By induction, we have all $T^{m}(u_i)$ are of good reduced forms
with denominators $u^{\gamma_{\Omega}(C^{m}(P_i[1]))}$ for any $m\ge
0, \ i=1,\cdots n.$ Dually, we have the same statement for
$T^{-m}(u_i)$, i.e. they are of good reduced forms with
denominators $u^{\gamma_{\Omega}(C^{-m}(P_i[1]))},$ where $m> 0.$
For any non-isomorphic objects $X, Y$ in $\mathcal{PI}(\Omega)$,
$\underline{\mbox{dim}}X\not=\underline{\mbox{dim}}Y$, then
$\phi_{\Omega}(X)\not=\phi_{\Omega}(Y).$  This says that
$\phi_{\Omega}$ is bijection.   The proof is finished.

\medskip

 If $(\G,\Omega)$ is a quiver (with trivial valuations), Caldero and Keller prove that the Caldero-Chapoton map $X_{?}$ [CC] gives
 a bijection from $\mathcal{E}(H)$  to
 $\chi_{B}$ in Theorem 4.7 [CK2]. In this case, our map $\phi_{\Omega}$ coincides with the map
$X_?$ since $X_?$ sends tilting objects to clusters, and then sends
$C^m(P_i[1])$ to $T^m(u_i).$ Therefore $\phi_{\Omega}$ sends cluster tilting
objects in $\mathcal{PI}$ to clusters in $\chi'_{\Omega}$.

\medskip

 In the following, we apply this theorem to the Dynkin diagrams  $\G$
 to obtain a generalization of Fomin and
Zelevinsky's denominator theorem [FZ2]. First of all we have a
direct consequence which is needed in the proof of the further
generalization (see the next proposition).
\medskip

\textbf{Corollary 4.6.} Let $\G$ be a Dynkin diagram and $M$ an
indecomposable $H-$module with $\underline{\mbox{dim}}M=\alpha$.
Then $\phi_{\Omega}(M)=\frac{f_{\alpha}(u_1,\cdots,
u_n)}{u^{\alpha}}, $ where $f_{\alpha}(u_1,\cdots, u_n)$ is an
integral polynomial and is not divided by $u_i$.

\medskip

\textbf{Remark 4.7.} From Corollary 4.6, combining with Lemma 8.2.
in [Ker], one can prove that different cluster monomials have
different denominators with respect to a given acyclic clusters (compare
[D]).

\medskip

Now we strengthen the denominator statement in Corollary 4.6 to any
 seed. Namely we generalize the denominator Theorem in [CCS2] to
non simply-laced Dynkin diagram by using the interpretation in [Z1]
of the compatibility degree $(\ ||\ ).$

\medskip

\textbf{Proposition 4.8.} Let $(\G,\Omega)$ be a valued Dynkin
quiver and $\underline{x}=(x_1,\cdots, x_n)$ a cluster of the
cluster algebra $\mathcal{A}_{\Omega}$. Let $V=\oplus_{i=1}^{i=n}V_i$
be the cluster tilting object corresponding to $\underline{x}$ under the map
$\phi_{\Omega}.$  Then  for any indecomposable object $M$ in
$\mathcal{C}(\Omega)$,
$\phi_{\Omega}(M)=\frac{f_{\alpha}(x_1,\cdots,
x_n)}{x^{\gamma_V(M)}}, $ where $f_{\alpha}(x_1,\cdots, x_n)$ is an
integral polynomial and is not divided by $x_i$.
\medskip

\textbf{Proof.} For simplicity, $\gamma_{\Omega}(M)$ is denoted by
$\alpha,$ $\gamma_{\Omega}(V_i[1])$ by $\beta_i.$ Let $-\alpha_i$ be
the negative simple roots. Since the bijection $\phi_{\Omega}$ sends
cluster tilting objects to clusters by Theorem 4.7. in [Z3], there is a cluster tilting
object $V$ which corresponds to the cluster $\underline{x}$ under
$\phi_{\Omega}.$ The cluster variable $\phi_{\Omega}(M)$ of cluster
algebra $\mathcal{A}_{\Omega}$ can be written as
$\frac{f_{\alpha}(x_1,\cdots, x_n)}{\prod_{i=1}^{i=n}x_i^{[\alpha,
\beta_i, \underline{x}]}}$ by the Laurent phenomenon [FZ2], where
$[\alpha, \beta_i, \underline{x}]\in \mathbf{Z}$. From Lemmas 6.2.
6.3. in [CCS1] (the proofs of these lemmas there were given for
simply-laced Dynkin case, but these proofs also work for non simply-laced Dynkin case without any changes), we have that $[\alpha,
\beta_i, \underline{x}]=[\alpha, \beta_i, \underline{y}]$ provided
$\beta_i\in \{x_j|j=1,\cdots n\}\bigcap \{y_j|j=1,\cdots n\}$ for
clusters $\underline{x}=(x_1,\cdots,x_n)$ and $\underline{y}=(y_1,\cdots, y_n)$. Then it is
denoted by $[\alpha, \beta_i]$ and $[\sigma _{\pm}\alpha, \sigma
_{\pm}\beta]=[\alpha,  \beta]$. From Corollary 4.6, we know that
$[\alpha, - \alpha_i]=n_i(\alpha)$ where $n_i(\alpha)$ is the $i-$th
coefficient of expression of root $\alpha $ in simple  roots. This
proves $[\alpha, - \alpha_i]=(-\alpha _i ||\alpha).$  By the
uniqueness of $(\ || \  )$ in Section 3 in [FZ2],  $[\alpha,
\beta]=( \beta || \alpha )$ for $\beta\not=\alpha.$  Then by Theorem
3.6 in [Z1], we have that $[\alpha,
\beta_i]=(\beta_i||\alpha)=\mbox{dim}_{\mbox{End}(V_i[1])}\mbox{Ext}^1_{\mathcal{C}(\Omega)}(V_i[1],
M)$ for $M\ncong V_i[1]$. The latter equals
$\mbox{dim}_{\mbox{End}(V_i[1])}\mbox{Hom}_{\mathcal{C}(\Omega)}(V_i[1],
M[1])=\mbox{dim}_{\mbox{End}(V_i[1])}\mbox{Hom}_{\mathcal{C}(\Omega)}(V_i,
M),$ it is the $i-$th component of dimension vector $\gamma_V(M).$
The proof is finished.

\medskip

\medskip

\begin{center}
\textbf {ACKNOWLEDGMENTS.}\end{center} This work was completed when
I was visiting Institut f\"ur Mathematik, Universit\"at Paderborn. I
thank Henning Krause and all members of  the group of representation
theory very much for warm hospitality and discussions. We would like
to thank Professor Idun Reiten for her helpful conservations on this
topic. The author is grateful to the referee for a number of helpful comments and valuable suggestions.

\begin{center}

\end{center}

\medskip

\begin{thebibliography}{99}

\bibitem[ABS]{abs}
I.Assem, T.~Br\"ustle and R.~Schiffler.
\newblock Cluster-tilted algebras as trivial extensions.
\newblock Preprint, {\tt arXiv:math.RT/0601537}, 2006.





\bibitem[BFZ]{bfz4}
A.Bernstein, S.~Fomin and A.~Zelevinsky.
\newblock Cluster algebras III: Upper bounds and double Bruhat cells.
\newblock Duke Math.J. \textbf{126}, no.1, 1-52, 2005.




\bibitem[BMR1]{ccs}
A.~Buan, R.Marsh, and I. Reiten.
\newblock Cluster-tilted algebras.
\newblock Preprint, {\tt arXiv:math.RT/0402054}, 2004. To appear in
Trans.AMS.


\bibitem[BMR2]{ccs}
A.~Buan, R.Marsh, and I. Reiten.
\newblock Cluster mutations via quiver representations.
\newblock Preprint, {\tt arXiv:math.RT/0412077}, 2004.


\bibitem[BMRT]{ccs}
A.~Buan, R.Marsh, I. Reiten and G.~Todorov.
\newblock Clusters and seeds in acycle cluster algebras.
\newblock Preprint, {\tt arXiv:math.RT/0510359}, 2005.


\bibitem[BMRRT]{ccs}
A.~Buan, R.Marsh, M.~Reineke,I. Reiten and G.~Todorov.
\newblock Tilting theory and cluster combinatorics.
\newblock  Advances in Math. \textbf{204}, 572-618, 2006.

\bibitem[CB]{}
W. Crawley-Boevey, Tame algebras and generic modules.
\newblock Proc.London Math. Soc., \textbf{63}, 241-264, 1991.

\bibitem[CC]{ccs}
P.~Caldero and F.~Chapoton.
\newblock Cluster algebras as Hall algebras of quiver representations.
\newblock Comment. Math. Helv. \textbf{81}, 595-616, 2006.

\bibitem[CCS1]{ccs}
P.~Caldero, F.~Chapoton and R.~Schiffler.
\newblock Quivers with relations arising from clusters ($A_n$ case).
\newblock \newblock Quivers with relations arising from clusters ($A_n$ case).
\newblock Transaction
of AMS. \textbf{358}, 1347-1364, 2006.



\bibitem[CCS2]{ccs}
P.~Caldero, F.~Chapoton and R.~Schiffler.
\newblock Quivers with relations and cluster tilted algebras.
\newblock Preprint {\tt arXiv:math.RT/0411238}, 2004.

\bibitem[CK1]{ck}
P.~Caldero and B.~Keller.
\newblock From triangulated categories to cluster algebras.
\newblock Preprint {\tt arXiv:math.RT/0506018}, 2005.


\bibitem[CK2]{ck}
P.~Caldero and B.~Keller.
\newblock From triangulated categories to cluster algebras II.
\newblock Preprint {\tt arXiv:math.RT/0510251}, 2005.



\bibitem[D]{d}
G.~Dupont.
\newblock An approach to non simply laced cluster algebras.
\newblock Preprint {\tt arXiv:math.RT/0512043 }, 2005.



\bibitem[FZ1]{fz1}
S.~Fomin and A.~Zelevinsky.
\newblock Cluster Algebras I: Foundations.
\newblock J. Amer. Math. Soc. \textbf{15}, no.2, 497--529, 2002.





\bibitem[FZ2]{fz3}
S.~Fomin and A.~Zelevinsky.
\newblock Cluster algebras II: Finite type classification.
\newblock Invent. Math. \textbf{154}, no.1, 63-121, 2003.



\bibitem[HRS]{}
D.~Happel, I.Reiten and S.Smal$\phi$.
\newblock Tilting in abelian categories and quasitilted algebras.
\newblock Mem. Amer. Math. Soc., \textbf{575}, 1996.




\bibitem[Ke]{ke}
B.~Keller.
\newblock Triangulated orbit categories.
\newblock Document Math. \textbf{10}, 551-581, 2005.


\bibitem[Ker]{ker}
O.~Kerner.
\newblock Representations of wild quiverss.
\newblock Canadian Mathematical Socieety Conference Proceedings Vol. 19, 65-107, 1996.


\bibitem[Kr]{kr}
H. Krause.
\newblock Stable equivalence preserves representation type.
\newblock Comment. Math.Helv. \textbf{72}, 266-284, 1997.


\bibitem[MRZ]{mrz}
R.~Marsh, M.~Reineke and A.~Zelevinsky.
\newblock Generalized associahedra via quiver representations.
\newblock Trans. Amer. Math. Soc. \textbf{355}, no.10, 4171-4186, 2003.







\bibitem[Z1]{z} B.Zhu.
\newblock BGP-reflection functors and Cluster combinatorics,
\newblock  Journal of Pure and Applied Algebra. In press. Also see {\tt arXiv:math.RT/0511380}

\bibitem[Z2]{z} B.Zhu.
\newblock Equivalences between cluster categories. \newblock  Journal of Algebra. In press.  Also see {\tt arXiv: math.RT/0511382}.

\bibitem[Z3]{z} B.Zhu.
\newblock Applications of BGP-reflection functors: isomorphisms of cluster algebras.
\newblock  Science in China. In press. Also see {\tt arXiv:math.RT/0511384}.


\end{thebibliography}
\end{document}